\documentclass[a4paper]{article}
\usepackage{amssymb, amsmath, amsthm, verbatim, cite, lscape, longtable, graphicx}
\usepackage[cp1251]{inputenc}
\usepackage[english]{babel}

\newtheorem{theorem}{Theorem}
\newtheorem{lemma}{Lemma}[section]
\newtheorem{prop}[lemma]{Proposition}
\newtheorem{cor}[lemma]{Corollary}

\theoremstyle{definition}
\newtheorem*{prf}{Proof}

\begin{document}
	
\begin{center}\textbf{\large A Criterion for the Existence of a Solvable $\pi$-Hall Subgroup in a Finite Group}
	
~
		
A.A. Buturlakin, A.P. Khramova
\end{center}

\textsc{Abstract}. Let $G$ be a finite group and let $\pi$ be a set of primes. In this paper, we prove a criterion for the existence of a solvable $\pi$-Hall subgroup of $G$, precisely, the group $G$ has a solvable $\pi$-Hall subgroup if, and only if, $G$ has a $\{p,q\}$-Hall subgroup for any pair $p$, $q\in\pi$.
\bigskip
	
\section{Introduction}

Let $\pi$ be a set of primes. Then $\pi'$ denotes a complement to $\pi$ in the set of all primes. A subgroup $H$ of a finite group $G$ is called a \textit{$\pi$-Hall subgroup} if its order $|H|$ is a \textit{$\pi$-number}, meaning all prime divisors of $|H|$ are in $\pi$, while the index $|G:H|$ is a $\pi'$-number. If a subgroup is $\pi$-Hall for some set of primes $\pi$, then it is a \textit{Hall subgroup}. 
	
In 1928, P. Hall has shown that a full counterpart of Sylow theorem for $\pi$-subgroups holds in any finite solvable group. It was proved later in 1937 by P.~Hall and independently in 1938 by S.A. \v{C}unihin that if a finite group $G$ has a $p'$-Hall subgroup for any prime $p$, then $G$ is solvable. In his paper \cite[p. 291]{Hall56} P. Hall has proposed a conjecture stating that a group which possesses a $\{p , q\}$-Hall subgroup for any prime $p$ and $q$ is solvable. Z. Arad and M. Ward have confirmed the conjecture in \cite{AradWard} with the use of the classification of finite simple groups. The result was later generalized by V.N.~Tyutyanov in 2002 \cite{Tyut02}, so that it is sufficient for a finite group to have a $\{2,p\}$-Hall subgroup for any prime $p$ in order to be solvable.
	
The primary goal of this paper is to prove the following generalization of Z.~Arad's and M.~Ward's result.
	
\begin{theorem}\label{t:main} Let $G$ be a finite group and $\pi$ be a set of primes. The group $G$ has a solvable $\pi$-Hall subgroup if, and only if, $G$ has a $\{p, q\}$-Hall subgroup for any~$p, q\in \pi$.
\end{theorem}

Observe that this statement cannot be strengthened in a fashion of Tyutyanov's result, that is the existence of $\{2,r\}$-Hall subgroups for every $r$ in $\pi$ does not guarantee the existence of a solvable $\pi$-Hall subgroup. Indeed, the projective special linear group $PSL_2(41)$ contains $\{2, 3\}$- and $\{2, 5\}$-Hall subgroups (see, \cite[Lemma 3.11]{RevVd10}), but does not contain ${\{2,3,5\}}$-Hall subgroup.
 
It is worth noting that \cite[Lemma 3.4]{Mor13} states that if $G$ is a finite group, $\pi$ is a set of primes containing at least three prime divisors of the order of $G$, and $G$ has a nilpotent $\tau$-Hall subgroup for any proper subset $\tau$ of $\pi$, then $G$ has a nilpotent $\pi$-Hall subgroup. This gives a weaker version of Theorem~\ref{t:main} for nilpotent Hall subgroups.

In order to prove Theorem~\ref{t:main} we first reduce it to the case of $G$ being an almost simple group. F.~Gross in \cite{Gross86} has obtained a sufficient condition for the existence of $\pi$-Hall subgroup in a finite group, precisely, such a subgroup exists if the group of induced automorphisms of any composition factor possesses a $\pi$-Hall subgroup of its own. However, it is not enough to assert the existence of a solvable $\pi$-Hall subgroup in a finite group.
	
The following result, partially analogous to F. Gross's \cite[Theorem 3.5]{Gross86}, hold for the solvable Hall subgroup case. The proof of the modified version is similar to the original for the most part.

\begin{theorem}\label{t:Reduction} Let $G$ be a finite group and $\pi$ be a set of primes. Let $$1=G_0\leqslant G_1\leqslant\dots\leqslant G_n=G$$ be a composition series of $G$ which is a refinement of a chief series of $G$. If group $Aut_G(G_i/G_{i-1})$ has a solvable $\pi$-Hall subgroup for any $1\leqslant i\leqslant n$, then $G$ has a solvable $\pi$-Hall subgroup.
\end{theorem}
	
The next step of the proof is a reduction to the simple group case. For a natural $n$ let $\pi(n)$ denote a set of its prime divisors, and for a finite group $G$ let $\pi(G)=\pi(|G|)$.
	
\begin{prop}\label{p:simple} Let $S$ be a finite simple group and $\pi$ be a set of primes such that $|\pi\cap\pi(S)|\geqslant 3$. The group $S$ has a $\{p, q\}$-Hall subgroup for any $p$, $q\in \pi$ if, and only if, the group $S$ has a solvable $\pi$-Hall subgroup. Furthermore, all the conjugacy classes of such subgroups are invariant with respect to the automorphism group.
\end{prop}
	
In view of Lemma~\ref{l:AlmSimpletoFinite} this proposition proofs Theorem~\ref{t:main} in case of $G$ being an almost simple group with $S$ as its socle, assuming $\pi\cap\pi(S)$ possesses at least three elements.
	
This paper is organized as follows. The first section describes preliminary results on Hall subgroups of finite groups. In the second section Theorem~\ref{t:Reduction} is proven. In sections three and four we prove Proposition~\ref{p:simple} and Theorem~\ref{t:main}.
	
\section{Notation and Preliminaries}
	
A subnormal subgroup $A$ of a group $G$ is an \textit{atom} if it is \textit{perfect}, i.e., it is equal to its commutator, and has a unique maximal normal subgroup denoted by $A^*$. By $\mathcal{A}(G)$ we denote the set of all atoms of a group $G$.
	
Suppose $H/K$ is a section of $G$. The normalizer $N_G(H/K)$ of this section is the intersection $N_G(H)\cap N_G(K)$, while its image in the automorphism group $Aut(H/K)$ is called a \textit{group of induced automorphisms} and is denoted by $Aut_G(H/K)$. If $K=1$, then we write $Aut_G(H)$ instead of $Aut_G(H/K)$.
	
\begin{lemma}\emph{\cite[Lemmas 2.1--2.3]{Gross86}} Let $H/K$ be a nonabelian composition factor of $G$.
		
$(1)$ There exists a unique $A\in\mathcal{A}(G)$ such that $H = AK$. Furthermore, $A^*=H\cap K$, $H/K\cong A/A^*$, $N_G(H/K)\leqslant N_G(A)$, and $Aut_G(H/K)$ isomorphic to a group $L$ with
$$Inn(A/A^*)\leqslant L \leqslant Aut_G(A/A^*).$$
		
$(2)$ $G^{(\infty)}= \langle A|A\in \mathcal{A}(G)\rangle.$
		
$(3)$ If $A\in\mathcal{A}(G)$ and $A\leqslant \langle H_1,\dots, H_m\rangle$, where $H_1,\dots, H_m$ are subnormal subgroups of $G$, then $A\leqslant H_i$ for some $i$, $1\leqslant i \leqslant m$.
\end{lemma}
	
Acquiring the notation introduced by P. Hall \cite{Hall56}, we denote the properties of $\pi$-subgroups of the group $G$ as follows.
	
$E_\pi$: the group $G$ has a $\pi$-Hall subgroup.
	
$C_\pi$: the group satisfies $E_\pi$ and all $\pi$-Hall subgroups are conjugate.
	
$D_\pi$: the group satisfies $C_\pi$, and every $\pi$-subgroup is contained in a $\pi$-Hall subgroup.
	
Equivalently to the wording ``$G$ satisfies $E_\pi$'', expressions like ``$G$ is an $E_\pi$-group'' or $G\in E_\pi$ may also be used. The same goes for the classes $C_\pi$ and $D_\pi$. In addition to that, $E_\pi^s$, $C_\pi^s$, and $D_\pi^s$ denote corresponding classes for solvable $\pi$-subgroups.
	
\begin{lemma} Suppose $M\trianglelefteq G$, $M\leqslant H\leqslant G$, $G\in E_\pi^s$, $H/M\in E_\pi$, and $G/M\in D_\pi$. Then $H$ has a solvable $\pi$-Hall subgroup. Furthermore, there exists a solvable $\pi$-Hall subgroup $A$ of $G$ such that $A\cap H$ is a solvable $\pi$-Hall subgroup of a group $H$.
\end{lemma}
	
\begin{prf} Follows from the proof of Lemma~3.1 in \cite{Gross86}.
\end{prf}
	
This statement has an immediate corollary.

\begin{cor}\label{c:Ginfty} Let $G\in E_\pi^s$ and $G^{(\infty)}\leqslant H\leqslant G$. Then $H\in E_\pi^s$.
\end{cor}
	
\begin{cor}\label{c:AutH(A/A*)} Assume $A\in \mathcal A(G)$, $A\leqslant H\leqslant G$, and $Aut_G(A/A^*)\in E_\pi^s$. Then $Aut_H(A/A^*)\in E_\pi^s$.
\end{cor}
	
\begin{prf} The proof follows one presented in~\cite[Corollary 3.3]{Gross86}. Since $Out(A/A^*)$ is solvable we have
$$Aut_G(A/A^*)^{(\infty)}\leqslant Inn(A/A^*)=Aut_A(A/A^*)\leqslant Aut_H(A/A^*).$$ Applying Corollary~\ref{c:Ginfty} to the group $Aut_G(A/A^*)$ completes the proof.
\end{prf}
	
The following two lemmas are of S. A. \v Cunihin.
	
\begin{lemma}\label{l:Cpi.Epi} Let $M\trianglelefteq G$. If $M\in C_\pi^s$ and $G/M\in E_\pi^s$, then $G\in E_\pi^s$.
\end{lemma}
	
\begin{prf} Without loss of generality we may assume $G/M$ to be a solvable $\pi$-group. Let $A$ be a solvable $\pi$-Hall subgroup of $G$. By Frattini argument $G=MN_G(A)$. The group $N_G(A)$ is $\pi$-solvable. Thus $N_G(A)\in D_\pi^s$.
\end{prf}
	
\begin{lemma}\label{l:AlmSimpletoFinite} Let $G$ be a finite group, $A\trianglelefteq G$, $\pi(G/A)\subseteq\pi$, and $M$ be a $\pi$-Hall subgroup of $A$. Then $G$ has a $\pi$-Hall subgroup $H$ such that $H\cap A=M$ if, and only if, $M^A=M^G$.
\end{lemma}
	
\begin{prf} Let us prove the ``if'' part. Since $(M^A)^G=M^A$, Frattini argument implies that $G=AN_G(M)$. The group $N_G(M)$ possesses a normal series
$$1\trianglelefteq M\trianglelefteq  N_A(M)\trianglelefteq N_G(M),$$
all sections of which are either $\pi$-groups or $\pi'$-groups. Thus, $N_G(M)$ is $\pi$-separable. It follows from~\cite[Theorem D6*]{Hall56} that $N_G(M)\in D_\pi$.
Therefore let $H$ be a $\pi$-Hall subgroup of $N_G(M)$ such that $M\leqslant H$.
From $\pi(G/A)\subseteq\pi$ and $\pi(|A:M|)\subseteq\pi'$ we derive that $|G:N_G(M)|$ is a $\pi'$-number. Hence $H$ is a $\pi$-Hall subgroup of $G$.
\end{prf}

\begin{lemma}\label{l:H.Hall(HM)} Assume $M\trianglelefteq G$, $G/M$ is a $\pi$-group, and $M=K_1\times\dots\times K_n$, where $K_1,\dots,K_n$ is a complete conjugacy class of subgroups in $G$. Let $N=N_G(K_1)$, $C=C_G(K_1)$ and $K=K_2\times\dots\times K_n$, and let $L$ be a subgroup of $N$ such that $L/C$ is a solvable $\pi$-Hall subgroup of~$N/C$.
Then there exists a $\pi$-Hall subgroup $H$ of $G$ such that $G=HM$, $H\cap S_1=L\cap S_1$, $L=(H\cap N)K$, and $H\cap M=(H\cap S_1)\times\dots\times(H\cap S_n)$.
\end{lemma}
\begin{prf}
Follows the proof of~\cite[Theorem 3.4]{Gross86}.
\end{prf}

The following results are used explicitly in the proof of Proposition~\ref{p:simple}. 

\begin{lemma}\emph{\cite[Theorem A4]{Hall56}}\label{l:Sym2primes} Let $S_n$ be a symmetric group and let $p<q\leqslant n$, where $p$ and $q$ are primes. Then $S_n\in E_{\{p,q\}}$ only when $p=2$, $q=3$, and $n\in\{3,4,5,7,8\}$.
\end{lemma}
	
\begin{lemma}\label{l:pinPi}\emph{\cite[Theorem 3.3]{Rev99} or \cite{RevVd11}} Let $G$ be a finite group of Lie type over a field of characteristic~$p$, and let $\pi$ be a set of primes such that~$p\in\pi$. If $H$ is a $\pi$-Hall subgroup of~$G$, then either $H$ is contained in a Borel subgroup, or it is a parabolic subgroup of~$G$.
\end{lemma}
	
\begin{lemma}\label{l:2pnotinPiClass}\emph{\cite[Theorem 4.9]{Gross95}} Let $G$ be a classical group over a field of characteristic~$p$, and let $\pi$ be a set of primes such that $2,p\notin\pi$. Then $G\in E_\pi$ if, and only if, $G\in E_{\{r,s\}}$ for each pair of distinct primes $r,s\in\pi$.
\end{lemma}
	
Let $q$ be a power of an odd prime $p$. Define $\varepsilon$ to be one of $\{+1,-1\}$ such that $q\equiv_4\varepsilon$ (henceforward $a\equiv b\pmod{n}$ is written as $a\equiv_n b$ for short). The letter $\varepsilon$ also denotes the corresponding sign of $\{+,-\}$. In similar fashion, the letter $\eta$ denotes a certain sign of the set $\{+,-,\circ\}$, whereby if $\eta\in\{+,-\}$ then $\eta$ also stands for the number $\eta1$.

\begin{lemma}\emph{\cite[Theorem 5.2]{RevVd06} or \cite[Theorem 8.9]{RevVd11}}\label{l:2inPi3pnotinPiClass} Let $G$ be a simple group of Lie type over a field of odd characteristic~$p$, and let $\pi$ be a set of primes such that $2\in\pi$, but $3,p\notin\pi$. Define $t\in\pi\setminus\{2\}$. Then $G\in E_\pi$ if, and only if, either $G=\,^2G_2(q)$ and $\pi=\{2,7\}$, or the following conditions hold:

$(1)$ $\pi\subseteq\pi(q-\varepsilon)$;

$(2)$ If $G$ is one of $PSL^\pm_n(q)$, $PSp^{}_{2n}(q)$, $P\Omega^\pm_{2n}(q)$, then $n<t$;
		
$(3)$ If $G=PSL^{-\varepsilon}_n(q)$, then $(n+1)/2<t$;
		
$(4)$ If $G=P\Omega^{-\varepsilon}_{2n}(q)$, then $n-1<t$ and $n$ is odd;
		
$(5)$ If $G=E^{-\varepsilon}_6(q)$, then $5\notin\pi$;
		
$(6)$ If $G=E_7(q)$ or $G=E_8(q)$, then $5,7\notin\pi$.
		
Furthermore, if $|\pi|\geqslant3$ then every $\pi$-Hall subgroup $H$ of a group $G$ has a normal $2$-complement, and all $\pi$-Hall subgroups are conjugate in~$G$.
\end{lemma}
	
Henceforward $n_\pi$ denotes $\pi$-part of a natural $n$, and $n_{\{p\}}$ is the same as $n_p$.
	
\begin{lemma}\label{l:E3sClass} Let $G$ be either $PSL_n^\pm(q)$ or $PSp_{2n}(q)$, where $q$ is a power of field characteristic $p$, and $s$ is a prime such that $s\notin\{2,p\}$, but $s\in\pi(q-\varepsilon)$. Then $G\in E_{\{3,s\}}$ if, and only if, one of the following holds:
		
$(1)$ $G=PSL_n(q)$, $e(q,3)=e(q,s)=a$, and $n<as$;
		
$(2)$ $G=PSL_3(q)$, $e(q,3)=2$, $e(q,s)=1$, and $(q^2-1)_3=3$;
		
$(3)$ $G=PSL^-_n(q)$, $e(q,3)=e(q,s)$, and $n<2s$;
		
$(4)$ $G=PSL^-_3(q)$, $e(q,3)=1$, $e(q,s)=2$, and $(q^2-1)_3=3$;
		
$(5)$ $G=PSp_{2n}(q)$, $e(q,3)=e(q,s)$, and $n<s$.
\end{lemma}

\begin{prf}
Follows immediately from~\cite[Theorems 4.1, 4.3, 4.5]{Gross95}.
\end{prf}
	
\section{Reduction to the Almost Simple Group Case}

This section, although not a direct corollary, mainly follows Gross's argument in \cite{Gross86}. Here we reprove his Theorem 3.4 and a certain part of Theorem 3.5 for them to hold for the $E_\pi^s$-groups instead of $E_\pi$-groups. We preserve most of Gross's argument.
Theorem~\ref{t:Reduction}, a modified version of Gross's Theorem 3.5, provides us with the desired reduction.

\begin{theorem}\label{t:CounterEx} Let $M\trianglelefteq G$ and assume that both $M$ and $G/M\in E_\pi^s$ but $G\notin E_\pi^s$. Then $G$ has a nonabelian composition factor $H/K$ with $H\leqslant M$ and $K\trianglelefteq G$ such that $H/K\in E_\pi^s$ but $Aut_G(H/K)\notin E_\pi^s$.
\end{theorem}

\begin{prf}
Assume we have a counterexample in which $|G|+|M|$ is as small as possible. Let $L$ be a normal subgroup of~$G$ such that $L<M$ and $M/L$ is a minimal normal in~$G/L$. Since $G\in E_\pi^s$, it follows that $L\in E_\pi^s$.
		
Assuming $G/L\in E_\pi^s$ we reach a contradiction to the minimality of the counterexample by replacing $M$ with~$L$. On the other hand, assuming $G/L\notin E_\pi^s$ we reach the same contradiction by factorizing $G$ and $M$ by~$L$. It follows that $L=1$, and $M$ is a chief factor of~$G$.
		
The group $M$ does not satisfy $C_\pi$. Because if it does, with $G/M\in E_\pi^s$ we have $G\in E_\pi^s$ by Lemma~\ref{l:Cpi.Epi}.

Therefore $M$ is not solvable. This implies that $M=K_1\times\dots\times K_n$, where $\{K_1,\dots,K_n\}$ is a complete conjugacy class of subgroups in~$G$, and $K_i$ for all~$i$ is a nonabelian composition factor of~$G$. Since $G\in E_\pi^s$ we have $K_i\in E_\pi$ for all~$i$.
		
In case $Aut_G(K_i)\notin E_\pi^s$, we could satisfy the conclusion of the theorem by choosing $K=1$ and $H=K_i$. Consider the case $Aut_G(K_i)\in E_\pi^s$ for all~$i$. It follows from Corollary~\ref{c:AutH(A/A*)} that $Aut_A(K_i)\in E_\pi^s$ for any~$A$ with $M\leqslant A\leqslant G$. It implies that $A\in E_\pi^s$ for any~$A$ with $M\leqslant A<G$ and $A/M\in E_\pi^s$, or else replacing $G$ with $A$ provides a counterexample to the theorem. However, $G/M\in E_\pi^s$. If $A/M$ is a $\pi$-Hall subgroup of~$G/M$, then $A\in E_\pi^s$ if, and only if, $G\in E_\pi^s$, which contradicts the hypothesis of the theorem. It follows that $G/M$ must be a $\pi$-group.
		
Now let $S=K_1$, $K=K_2\times\dots\times K_n$, $N=N_G(S)$, and $C=C_G(S)$.
Let $L/C$ be a solvable $\pi$-Hall subgroup of $N/C\in E_\pi^s$. Since $G/M$ is a $\pi$-group, we must have $N=LM$.
		
We now use Lemma~\ref{l:H.Hall(HM)} to claim that~$G$ contains a $\pi$-Hall subgroup $B$ such that $G=BM$, $L=(B\cap N)K$, $B\cap S=L\cap S$, and $B\cap M=(B\cap K_1)\times\dots\times(B\cap K_n)$.
		
Consider $(L\cap S)C/C$, a subgroup in~$L/C$. It is isomorphic to $(L\cap S)/(L\cap S\cap C)$, but $S\cap C=1$. Then $B\cap S = L\cap S$ is a solvable $\pi$-Hall subgroup. Since $G/M=BM/M\in E_\pi^s$, the group $B$ is solvable, and $G\in E_\pi^s$, contrary to the hypothesis.
\end{prf}
	
Now we prove Theorem~\ref{t:Reduction}.
	
\begin{prf}
By~\cite[Lemma 2.5]{Gross86} we have $Aut_G(A/A^*)\in E_\pi^s$ for all $A\in\mathcal{A}(G)$.

We prove by induction on~$|G|$. Let $M$ be a minimal normal subgroup of~$G$. From~\cite[Lemma 2.4]{Gross86} we find that $Aut_{G/M}(A/A^*)\in E_\pi^s$ where $A\in\mathcal{A}(G/M)$. By induction, $G/M\in E_\pi^s$. If $M$ is abelian, then by Lemma~\ref{l:Cpi.Epi} we obtain $G\in E_\pi^s$. If $M$ is not abelian, then $M=A_1\times\dots\times A_n$, where each $A_i$ is a nonabelian simple group. Then $A_i\in\mathcal{A}(G)$ for all~$i$.
Since by assumption $Aut_G(A_i)\in E_\pi^s$, and since $Inn(A_i)\trianglelefteq Aut_G(A_i)$, we see that $A_i\in E_\pi^s$ and hence $M\in E_\pi^s$.
By Theorem~\ref{t:Reduction} we have $G\in E_\pi^s$, and the proof is complete.
\end{prf}
	
\section{The Proof of Proposition~\ref{p:simple}}
	
Let $S$ be either symmetric or alternating group. Then by Lemma~\ref{l:Sym2primes} the group $S$ does not have any $\{p,q\}$-Hall subgroups unless $p=2$ and $q=3$.
	
Let $S$ be a sporadic group. Then $S\in E_{\{p,q\}}$ for any $p,q\in\pi$ only if either $\pi=\{p,q\}$, or $G=J_1$ and $\pi=\{2,3,7\}$. Either way, $S\in E_\pi^s$.
	
Let $S$ be a simple group of Lie type over a field of characteristic $p$. Let us first consider the case of $p\in\pi$.
	
\begin{prop}\label{p:pinPi} Let $S$ be a finite simple group of Lie type over a field of characteristic~$p$, and suppose $B$ is a Borel subgroup of~$S$, and $\pi$ is a set of primes such that $p\in\pi$. If $|\pi\cap\pi(S)|\geqslant3$ and~$S$ has $\{r, s\}$-Hall subgroup for any $r, s\in\pi$, then $|S|_\pi=|B|_\pi$.
\end{prop}
	
\begin{prf} Define $\pi^*=\pi\cap\pi(S)\setminus\{p\}$. According to Lemma~\ref{l:pinPi} we only need to show that no $\{r, s\}$-Hall subgroup of~$S$ cannot be parabolic.
		
Assume, on the contrary, that $\{p, r\}$-Hall subgroup is parabolic for some $r\in\pi^*$. Then $\{p, s\}$-Hall subgroup is parabolic for any $s\in\pi^*$. Indeed, by Lemma~\ref{l:pinPi} if $\{p,s\}$-Hall subgroup is not parabolic, then it is contained in Borel subgroup, and therefore in any parabolic subgroup. In particular, it should be contained in $\{p, r\}$-Hall parabolic subgroup, which is not possible.
		
However, the order of a parabolic subgroup which is not a Borel subgroup is divisible by $6$; the contradiction is reached.
\end{prf}
	
Thus, if $S\in E_{\{r,s\}}$ for any $r,s\in\pi$, then $\pi$-Hall subgroup $H$ of $S$ is in a Borel subgroup $B$, therefore $H$ must be solvable. Moreover, since any $\pi$-Hall subgroup is in $B$, all $\pi$-Hall subgroups are conjugate.
	
The following lemmas cover the case when characteristic~$p$ is not in~$\pi$.
	
\begin{lemma}\label{l:2pnotinPi} Let $G$ be a simple group of Lie type over a field of characteristic~$p$, and let $\pi$ be a set of primes such that $2,p\notin\pi$. Then $G\in E_\pi$ if, and only if, $G\in E_{\{r,s\}}$ for any $r,s\in\pi$.
\end{lemma}
	
\begin{prf}	Suppose $G$ is classical group. Then the statement holds by Lemma~\ref{l:2pnotinPiClass}.
		
Let $G$ be a Rie or Suzuki group. Then according to~\cite[Lemma 14]{RevVd02} the group $G$ satisfies $E_\pi$ if $\pi\cap\pi(G)$ is contained in one of the sets listed in the table:
\begin{center} Table
			
\begin{tabular}{|c|c|}\hline
$G$& $\pi\cap\pi(G)\subseteq$\\ \hline
$^2B_2(2^{2n+1})$ & $\pi(2^{2n+1}-1)$\\
& $\pi(2^{2n+1}\pm2^{n+1}+1)$ \\ \hline
$^2G_2(3^{2n+1})$ & $\pi(3^{2n+1}-1)$ \\
& $\pi(3^{2n+1}\pm3^{n+1}+1)$ \\ \hline
$^2F_4(2^{2n+1})$ & $\pi(2^{2(2n+1)}\pm1)$ \\
& $\pi(2^{2n+1}\pm2^{n+1}+1)$ \\
& $\pi(2^{2(2n+1)}\pm2^{3n+2}\mp2^{n+1}-1)$ \\
& $\pi(2^{2(2n+1)}\pm2^{3n+2}+2^{2n+1}\pm 2^{n+1}-1)$\\ \hline
\end{tabular}
\end{center}

By a straightforward calculation it can be verified that corresponding sets of prime divisors intersect trivially. Thus, the ''if'' part of the statement is true.
		
If $G$ is one of $^3D_4(q)$, $E^\pm_6(q)$, $E_7(q)$, $E_8(q)$, $F_4(q)$, $G_2(q)$, then according to~\cite[Lemmas 7--13]{RevVd02} $G$ satisfies $E_\pi$ if, and only if, the following conditions hold:
		
$(1)$ For $r=\min\pi\cap\pi(G)$ and for any $s\in\pi\cap\pi(G)\setminus\{r\}$ we have $e(q,r)=e(q,s)$;
		
$(2)$ If $G=E^\pm_6(q)$, then $(q\mp1)_\pi\not\equiv_{15}0$;
		
$(3)$ If $G=E_7(q)$, then $e(q,r)=1$ and $(q-1)_\pi$ is not divisible by $15$, $21$, or $35$;
		
$(4)$ If $G=E_8(q)$, then $e(q,r)=2$ and $(q+1)_\pi$ is not divisible by $15$, $21$, or $35$.
		
It is evident that the first condition holds for any pair $r,s\in\pi$ if, and only if, it holds for $\pi$.
		
Let $G=E_7(q)$. The conditions $(q-1)_{\{r,s\}}\not\equiv_{15}0$, $(q-1)_{\{r,s\}}\not\equiv_{21}0$, and $(q-1)_{\{r,s\}}\not\equiv_{35}0$ hold for any $r,s\in\pi$ if, and only if, $\pi$ has at most one element of the set $\{3,5,7\}$. Then $(q-1)_\pi\not\equiv_k0$, where $k\in\{15,21,35\}$.

In cases of $G=E^\pm_6(q)$ and $G=E_8(q)$ the argument is analogous to one in previous paragraph.
\end{prf}
	
\begin{lemma}\label{l:2inPi3pnotinPi} Let $G$ be a simple group of Lie type over a field of characteristic~$p$, and let $\pi$ be a set of primes such that $2\in\pi$ and $3,p\notin\pi$. Then $G\in E_\pi^s$ if, and only if, $G\in E_{\{r,s\}}$ for any $r,s\in\pi$. Furthermore, all solvable $\pi$-Hall subgroups are conjugate in~$G$.
\end{lemma}
	
\begin{prf}
By Lemma~\ref{l:2inPi3pnotinPiClass}, $G\in E_{\{r,s\}}$ for any $r,s\in\pi$ if, and only if, $G\in E_\pi$. Furthermore, if $|\pi|\geq3$ then every $\pi$-Hall subgroup $H$ of $G$ has a normal $2$-complement~$K$. Since a subgroup of odd order is solvable, and $H/K$ is a cyclic $2$-group, we can conclude that $H$ is solvable. In addition, whenever $|\pi|\geqslant3$ all $\pi$-Hall subgroups are conjugate in~$G$ by the same Lemma~\ref{l:2inPi3pnotinPiClass}.
\end{prf}
	
\begin{lemma}\label{l:23inPipnotinPi} Let $G$ be a simple group of Lie type over a field of characteristic~$p$, and let $\pi$ be a set of primes such that $2,3\in\pi$ and $p\notin\pi$. Then $G\in E_\pi^s$ if, and only if, $G\in E_{\{r,s\}}$ for any $r,s\in\pi$. Furthermore, if $|\pi\cap\pi(G)|>2$, then the conjugacy classes of solvable $\pi$-Hall subgroups are invariant with respect to the automorphism group.
\end{lemma}
	
\begin{prf} Define $\tau=\pi\setminus\{ 2, 3\}$. Henceforward $s$ denotes some element of~$\tau$.

Let $G=PSL^\eta_2(q)$. Since $G\in E_{\{2,s\}}$, it follows from Lemma~\ref{l:2inPi3pnotinPiClass} that $s$ is in $\pi(q-\varepsilon)$. If $3$ divides $q-\varepsilon$ then by~\cite[Lemma 3.11]{RevVd10} the group $G$ has a solvable $\pi$-Hall subgroup and all such subgroups are conjugate. Suppose $3$ does not divide $q-\varepsilon$ and $G$ has a $\{ 3, s\}$-Hall subgroup. According to Lemma~\ref{l:E3sClass}, if group $PSL_2^\eta(q)$ has a $\{ 3, s\}$-Hall subgroup, then $e(q ,3)= e(q, s)$, which, with regard for $s\in\pi(q-\varepsilon)$, contradicts the suggestion of $3$ not dividing~$q-\varepsilon$.

Let $n\not=2$. By~\cite[Lemma 4.3]{RevVd10}, to prove that $G\in E^s_\pi$ it is sufficient to show that one of these two conditions holds:
		
\begin{center}
	\begin{tabular}{lc}
	$(\text{A})$ & $q\equiv_{12}\eta \text{ or } n=3 \text{ and } q\equiv_4\eta;$ \\
	& $S_n\in E_\pi;$ \\
	& $\pi\cap\pi(G)\subseteq\pi(q-\eta)\cup\pi(n!);$ \\
	& $\text{if } r\in\pi\cap\pi(n!)\setminus\pi(q-\eta), \text{ then } |G|_r=|S_n|_r.$\\
	&\\
	$(\text{B})$ & $q\equiv_3-\eta;$ \\
	& $S_m\in E_\pi,\text{ where } n=2m\text{ or }n=2m+1;$ \\
	& $GL^\eta_2(q)\in E_\pi;$ \\
	& $\pi\cap\pi(G)\subseteq\pi(q^2-1).$\\
\end{tabular}	
\end{center}	

Since $G\in E_{\{2,s\}}$, it follows from Lemma~\ref{l:2inPi3pnotinPiClass} that $s\in\pi(q-\varepsilon)$, and either $n<s$ or $\eta=-\varepsilon$ and $\frac{n+1}{2}<s$.
Since $G\in E_{\{2,3\}}$, by~\cite[Lemma 4.3]{RevVd10} one of the conditions $(1-3)$ holds.
		
$(1)$ $q\equiv_{12}\eta\text{ or }n=3\text{ and }q\equiv_4\eta; S_n\in E_{\{2,3\}};\text{ if }3\notin\pi(q-\eta),\text{ then }|G|_3=3.$
		
In particular, it follows from $(1)$ that $\eta=\varepsilon$, therefore $n<s$. Then $S_n\in E_\pi$. Since $s\in\pi(q-\varepsilon)$ for any $s\in\tau$, we have $(\pi\cap\pi(n!))\setminus\pi(q-\eta)\subseteq\{3\}$. In the last formula the equality holds only if $n=3$, and in this case $|G|_3=|S_3|_3=3$. So, the condition $(\text{A})$ holds and $G\in E_\pi^s$.
		
$(2)$ $q\equiv_3-\eta; S_m\in E_{\{2,3\}},\text{ where } n=2m\text{ or }n=2m+1; GL^\eta_2(q)\in E_{\{2,3\}}.$
		
By~\cite[Lemma 3.2]{RevVd10}, it follows from $GL^\eta_2(q)\in E_{\{2,3\}}$ that $3\in\pi(q-\varepsilon)$ or $(q^2-1)_{2,3}=24$.
		
Let $3\in\pi(q-\varepsilon)$. Since $s\in\pi(q-\varepsilon)$, by~\cite[Lemma 3.2]{RevVd10} we have $GL^\eta_2(q)\in E_\pi$.
If $\eta=\varepsilon$, then $m<n<s$ and $S_m\in E_\pi$. Let $\eta=-\varepsilon$. Since $\{3,s\}\subseteq\pi(q-\varepsilon)$, then $e(q,3)=e(q,s)$ for any $s\in\tau$, and since $G\in E_{\{3,s\}}$, $n<2s$ by Lemma~\ref{l:E3sClass}. Then $m<s$, thus $S_m\in E_\pi$.
So, $q\equiv_3-\eta$, $S_m\in E_\pi$, $GL^\eta_2(q)\in E_\pi$ and $\pi\cap\pi(G)\subseteq\pi(q^2-1)$ since $s\in\pi(q-\varepsilon)$. Thus, $(\text{B})$ holds and $G\in E_\pi^s$.

Let $3\notin\pi(q-\varepsilon)$ and $(q^2-1)_{2,3}=24$. From the fact that $3$ divides $(q+\eta)$ we obtain that $\eta=\varepsilon$. In particular, $e(q,3)\not=e(q,s)$ and from $G\in E_{\{3,s\}}$ by Lemma~\ref{l:E3sClass} we have $n=3$. Therefore $\pi\cap\pi(G)\subseteq\pi(q-\eta)\cup\pi(n!)$ and $q\equiv_4\eta$. Also, if $n=3$ and $(q^2-1)_3=3$ then $|G|_3=3$, thus, for any $r\in(\pi\cap\pi(n!))\setminus\pi(q-\eta)=\{3\}$ we have $|G|_r=|S_3|_r=3$. Finally, since $S_3\in E_\pi$, the condition $(\text{B})$ holds and $G\in E_\pi^s$.
		
$(3)$ $n=11; (q^2-1)_{2,3}=24; q\equiv_3-\eta; q\equiv_4\eta.$
		
If $3\notin\pi(q-\varepsilon)$, then $e(q,3)\not=e(q,s)$, and since $G\in E_{\{3,s\}}$, by Lemma~\ref{l:E3sClass} we have $n=3$. Thus, $3\in\pi(q-\varepsilon)$. Then $\pi\subseteq\pi(q-\varepsilon)$ and by~\cite[Lemma 3.2]{RevVd10} it follows that $GL^\eta_2(q)\in E_\pi$. Since $q\equiv_3-\eta$, $\pi\cap\pi(G)\subseteq\pi(q^2-1)$ and $S_5\in E_\pi$, the group $G$ has solvable $\pi$-Hall by the condition~$(\text{B})$.
		
Let $G=PSp_{2n}(q)$. By~\cite[Lemma 4.4]{RevVd10}, to prove that $G$ is an $E^s_\pi$-group it suffices to show that $S_n\in E_\pi$, $SL_2(q)\in E_\pi$, and $\pi\cap\pi(G)\subseteq\pi(q^2-1)$.
Since $G\in E_{\{2,s\}}$, by Lemma~\ref{l:2inPi3pnotinPiClass} we have $s\in\pi(q-\varepsilon)$ and $n<s$.
Since in addition to that $G\in E_{\{3,s\}}$, by Lemma~\ref{l:E3sClass} the equality $e(q,3)=e(q,s)$ holds, thus $3\in\pi(q-\varepsilon)$, and then by~\cite[Lemma 3.1]{RevVd10} we have $SL_2(q)\in E_\pi$.
Since $G\in E_{\{2,3\}}$, it follows from~\cite[Lemma~4.4]{RevVd10} that $S_n\in E_{\{2,3\}}$, and since $n<s$ we obtain $S_n\in E_\pi$.
Finally, since the condition $\pi\cap\pi(G)\subseteq\pi(q^2-1)$ holds as well, $G$~has a solvable $\pi$-Hall subgroup.
		
Let $G=P\Omega^\eta_n(q)$, and $m$ is such that $n=2m$ or $n=2m+1$. By~\cite[Lemma 6.7]{RevVd10}, to prove the existence of solvable $\pi$-Hall subgroup it suffices to show that one of these three conditions holds.
\begin{center}
\begin{tabular}{lc}
$(\text{A})$ & $n\text{ is odd}; \pi\cap\pi(G)\subseteq\pi(q-\varepsilon); q\equiv_{12}\varepsilon; S_m\in E_\pi.$ \\

$(\text{B})$ & $n\text{ is even}; \eta=\varepsilon^m; \pi\cap\pi(G)\subseteq\pi(q-\varepsilon); q\equiv_{12}\varepsilon; S_m\in E_\pi.$ \\

$(\text{C})$ & $n\text{ is even}; \eta=-\varepsilon^m; \pi\cap\pi(G)\subseteq\pi(q-\varepsilon); q\equiv_{12}\varepsilon; S_{m-1}\in E_\pi.$ \\
\end{tabular}
\end{center}

Since $G\in E_{\{2,3\}}$, $m\geqslant3$ by~\cite[Lemma 6.7]{RevVd10}. From $G\in E_{\{2,s\}}$ with the use of Lemma~\ref{l:2inPi3pnotinPiClass} we obtain that $s\in\pi(q-\varepsilon)$ and $m<s$, or if $\eta=-\varepsilon$ then $m$~is odd and $m-1<s$.
		
By~\cite[Lemma 6.7]{RevVd10}, $G\in E_{\{2,3\}}$ if one of the conditions $(1-5)$ holds.
		
$(1)$ $n=2m+1$; $3\in\pi(q-\varepsilon)$; $S_m\in E_{\{2,3\}}$.
		
Then $\pi\cap\pi(G)\subseteq\pi(q-\varepsilon)$, and $S_m\in E_\pi$ since $m<s$. Thus, the condition~$(\text{A})$ holds and $G\in E_\pi^s$.
		
$(2)$ $n=2m$; $\eta=\varepsilon^m$; $3\in\pi(q-\varepsilon)$; $S_m\in E_{\{2,3\}}$.
		
If $\eta=-\varepsilon$, then $m$ is odd, which contradicts $\eta=\varepsilon^m$. We derive that $\eta=\varepsilon$, so $m<s$ and $S_m\in E_\pi$. Thus, $(\text{B})$ holds and $G\in E_\pi^s$.
		
$(3)$ $n=2m$; $\eta=-\varepsilon^m$; $3\in\pi(q-\varepsilon)$; $S_{m-1}\in E_{\{2,3\}}$.
		
Then $m-1<s$, so $S_{m-1}\in E_\pi$, and in sight of $(\text{C})$ the group $G$ has a solvable $\pi$-Hall subgroup.
		
$(4)$ $n=11$; $3\in\pi(q-\varepsilon)$; $(q^2-1)_{2,3}=24$.
		
Since $\pi\cap\pi(G)\subseteq\pi(q-\varepsilon)$ and $S_5\in E_\pi$, the condition $(\text{A})$ holds and $G\in E_\pi^s$.
		
$(5)$ $n=12$; $\eta=-1$; $3\in\pi(q-\varepsilon)$; $(q^2-1)_{2,3}=24$.
		
Since $m=6$ is even, $\eta=\varepsilon=-1$, and $\eta=-\varepsilon^m$. Also, $\pi\cap\pi(G)\subseteq\pi(q-\varepsilon)$ and $S_5\in E_\pi$, therefore $(\text{C})$ holds, and $G$ has a solvable $\pi$-Hall subgroup.
		
Let $G$ be an exceptional group. By~\cite[Lemmas 7.1--7.6]{RevVd10}, to prove $G\in E^s_\pi$ it is sufficient to show that $\pi\cap\pi(G)\subseteq\pi(q-\varepsilon)$. Since $G\in E_{\{2,3\}}$, by~\cite[Lemmas 7.1--7.6]{RevVd10} we obtain that $3\in\pi(q-\varepsilon)$, and additionally $G$ has to be one of the groups $F_4(q)$, $G_2(q)$, $^3D_4(q)$, or $E^{-\varepsilon}_6(q)$. By Lemma~\ref{l:2inPi3pnotinPiClass}, all these groups satisfy $E_{\{2,s\}}$ if $s\in\pi(q-\varepsilon)$. Then $\pi\cap\pi(G)\subseteq\pi(q-\varepsilon)$, and the proof is complete.
\end{prf}

\section{The Proof of Theorem~\ref{t:main}}
	
Let $G\in E_{\{p,q\}}$. Then, by the $E_\pi$ criterion~\cite[Corollary 5]{RevVd11.Epi}, we have $Aut_G(G_i/G_{i-1})$ is in $E_{\{p,q\}}$ for any~$i$.
On the other hand, by Theorem~\ref{t:Reduction}, if $Aut_G(G_i/G_{i-1})\in E_\pi^s$, then $G\in E_\pi^s$.
Thus, to prove the main theorem it is sufficient to consider the case of $G$ being an almost simple group.
	
Let $G$ be an almost simple group, whereas $S=Soc(G)$.
Let $\pi\cap\pi(S)=\{r,s\}$ and $\pi\cap\pi(G)=\{r,s\}\cup\tau$. Since $G\in E_{\{r,s\}}$, there exists a $\pi$-Hall subgroup $H$ in~$S$ such that $N_G(H)\in E_{\{r,s\}}$. Let $T$ be a $\tau$-Hall subgroup of~$G$. Since the orders of $T$ and $S$ are coprime, by Glauberman's lemma $T\leqslant N_G(H)$. Thus, the $\pi$-Hall subgroup of $G$ is also $\pi$-Hall subgroup in $N_G(H)$. Then $|H^G|=|G:N_G(H)|$ is a $\pi'$-number, and $|H^S|=|S:N_S(H)|$, therefore $G/S\cong N_G(H)/N_G(S)$ and $H^G=H^S$, so now Lemma~\ref{l:AlmSimpletoFinite} can be used.
	
Let $|\pi\cap\pi(S)|>2$. By Proposition~\ref{p:simple}, the statement of the theorem holds for~$S$, and all the conjugacy classes of solvable $\pi$-Hall subgroups are invariant with respect to~$G$. Thus, by Lemma~\ref{l:AlmSimpletoFinite} the statement of the theorem is true for~$G$ as well.


\begin{thebibliography}{}

\bibitem{AradWard} Z.~Arad, M.~Ward. New criteria of solvability of finite groups. {\it J. Algebra} \textbf{77} (1982), 234--246.
\bibitem{Gross86} F.~Gross. On the Existence of Hall subgroups. {\it J. Algebra} \textbf{98}:1 (1986), 1--13.
\bibitem{Gross95} F.~Gross. Odd order Hall subgroups of the classical linear groups. {\it Math. Z.} 220 (1995), 317--336.
\bibitem{Hall56} P.~Hall. Theorems like Sylow’s. {\it Proc. London Math. Soc.} (3) \textbf{6} (1956), 286--304.
\bibitem{Mor13} A.~Moret\'{o}. Sylow numbers and nilpotent Hall subgroups. {\it J. Algebra} \textbf{379} (2013), 80--84.
\bibitem{Rev99} D.\,O.~Revin. Hall $\pi$-subgroups of finite Chevalley groups whose characteristic belongs to $\pi$. {\it Sib. Adv. Math.} 9 No. 2 (1999), 25-71.
\bibitem{RevVd06} D.\,O.~Revin, E.\,P.~Vdovin. Hall subgroups of finite groups. {\it Ischia group theory 2004: Proceedings of a conference in honor of Marcel Herzog (Naples, Italy, 2004), Contemp. Math.} \textbf{402}, {\it Amer. Math. Soc., Providence, RI; Bar-Ilan University, Ramat Gan} (2006), 229--263.
\bibitem{RevVd10} D.\,O.~Revin, E.\,P.~Vdovin. On the number of classes of conjugate Hall subgroups in finite simple groups. {\it J. Algebra} \textbf{324}:12 (2010), 3614--3652.
\bibitem{RevVd11.Epi} D.\,O.~Revin, E.\,P.~Vdovin. Existence criterion for Hall subgroups of finite groups. {\it J. Group Theory} \textbf{14}:1 (2011), 93--101.
\bibitem{RevVd02} E.\,P.~Vdovin, D.\,O.~Revin. Hall subgroups of odd order in finite groups. {\it Algebra and Logic} \textbf{41}:1 (2002), 8--29.
\bibitem{RevVd11}E.\,P.~Vdovin, D.\,O.~Revin. Theorems of Sylow type. {\it Russian Math. Surveys} \textbf{66}:5 (2011), 829--870.
\bibitem{Tyut02}V.\,N.~Tyutyanov. On a Hall hypothesis. {\it Ukrainian Math. J.} \textbf{54}:7 (2002), 1181--1191.
\end{thebibliography}
\end{document}